# Data-driven identification of critical links in transport networks using quantum annealing

Junxiang Xu[1], Chence Niu[2,*], Tingting Zhang[3,4], Divya Jayakumar Nair[1], Vinayak Dixit[1]

1. Research Centre for Integrated Transport Innovation (rCITI), School of Civil and Environmental Engineering, The University of New South Wales, Kensington, UNSW Sydney, NSW, 2052, Australia
2. Guangdong Basic Research Center of Excellence for Ecological Security and Green Development, Key Laboratory for City Cluster Environmental Safety and Green Development of the Ministry of Education, School of Ecology, Environment and Resources, Guangdong University of Technology, Guangzhou, 510006, China
3. Key Laboratory for Vehicle Emission Control and Simulation of Ministry of Ecology and Environment, Chinese Research Academy of Environmental Sciences, Beijing, China
4. Vehicle Emission Control Center, Chinese Research Academy of Environmental Sciences, Beijing, China

**corresponding author*:** niuchence@163.com

**Abstract:** In urban transport systems, time-varying demand and network conditions cause the importance of infrastructure elements to evolve, requiring the identification of period-specific critical links to support system-level risk and resilience analysis. However, static or time-averaged network analyses struggle to capture the temporal variation of infrastructure importance at the city scale. To address this gap, this study proposes a time-dependent critical link identification framework for large-scale urban transport networks. The problem is formulated as a Quadratic Unconstrained Binary Optimisation (QUBO) model and solved using quantum annealing on D-Wave hardware. Empirical analysis using real-world traffic data reveals a strong temporal concentration of critical links. Rather than persistently influencing system performance, critical links emerge mainly within a small number of key time windows, during which even limited disruptions can lead to substantial network delay amplification. These findings demonstrate the value of time-dependent analysis for risk screening, stress testing, and resilience-oriented transport management.

## 1. Introduction

As urban systems continue to expand in scale, the high level of interconnectivity among urban infrastructure components means that local failures can rapidly escalate into system wide risks (Zimmerman, 2001; Little, 2002). Against this background, accurately identifying critical infrastructure elements that play a decisive role in urban operations has become a central issue in urban resilience planning, risk assessment, and emergency management (Monstadt and Schmidt, 2019). Critical urban infrastructure generally refers to network nodes or links whose failure may trigger substantial functional degradation, cascading failures, or service disruption (Bodenstein and Pfeffer, 2026). The identification of such elements directly informs the prioritisation of limited resources, the formulation of risk mitigation strategies, and the overall stability of urban systems under extreme scenarios (Gu et al., 2023). A key question therefore concerns how to systematically identify critical infrastructure elements within complex urban networks, which carries important theoretical significance and practical value for enhancing urban capacity to cope with uncertain shocks.

*1.1 Research background and motivation*

Urban transport networks cover extensive spatial areas, contain many links, and operate under conditions that evolve continuously over time (Lin and Ban, 2013). Congestion, incidents, construction activities, or local infrastructure failures at any given moment can propagate through route substitution and flow redistribution, ultimately affecting the continuity of urban functions and the reliability of travel services (Saberi et al., 2020). In such systems, overall network performance is often determined not by the average contribution of all links, but by a small set of links that serve as critical corridors during specific time periods (Almotahari and Yazici, 2021). Therefore, critical link identification for urban transport infrastructure should not focus on producing a static list of important links. Instead, it should capture how criticality changes across different traffic states while remaining computationally tractable and practically applicable for large-scale networks (Huang et al., 2023; Alizadeh and Dodge, 2025).

The core motivation of this study is to develop a framework that can rapidly identify critical links in large-scale urban transport networks and characterise the temporal variation of their importance under time-dependent conditions. By extending critical link identification from a single time or a time averaged perspective to the full time series, this study seeks to more clearly reveal which links are most sensitive to network delay and performance degradation, and at which times. This perspective supports more targeted time sensitive monitoring, maintenance planning, and risk response decisions.

*1.2 Current research status and gaps*

Research on the identification of critical elements in urban infrastructure has mainly focused on two strands, namely network structural analysis and operational performance evaluation. One strand examines topological characteristics and measures structural importance using indicators such as degree, betweenness, and connectivity (Yu et al., 2017; Yang et al., 2020; Liu et al., 2024), recent topology-structure studies advance critical link identification by integrating travel demand under network congestion scenarios (Almotahari and Yazici, 2020; Sugiura and Kurauchi, 2023; Nitheesh and Bhavathrathan, 2025). The other strand incorporates traffic assignment or operating state information and evaluates the impact of facility failures on network performance using measures such as delay, capacity loss, or changes in system efficiency (Ben-Akiva et al., 2012; Shafiei et al., 2018; Long and Szeto, 2019). Such methods have shown strong recent application in urban road network design and post-disaster road recovery research (Pei et al., 2024; Wang et al., 2026). In summary, these current approaches provide an important foundation for understanding the role of critical elements in urban transport networks.

However, in real world urban operating environments, network criticality is not a static property (Cats and Jenelius, 2014). It varies substantially over time. Traffic demand levels, route choice behaviour, and network load distribution differ across time periods, which means that the impact of the same link on system performance can change markedly from one moment to another (de Moraes Ramos et al., 2020). Analyses that rely only on a single time or on time averaged conditions often struggle to accurately capture the true importance of critical links under high load or abnormal operating states (Li et al., 2025).

At the same time, as urban transport networks continue to expand in scale, critical link identification faces increasing computational challenges (Lin and Ban, 2013; Zhou and Wang, 2017; Huang et al., 2023; Wang et al., 2025). Systematically accounting for multi-link disruption combinations and their time-dependent characteristics in large-scale networks often requires trade-offs between modelling accuracy and computational efficiency (Almotahari and Yazici, 2021; Munikoti et al., 2022). This constraint has limited the application of existing methods across full temporal horizons and city scale settings. This constraint has limited application of existing methods across full temporal horizons and city scale settings. How critical links and infrastructure influence transport network resilience under disaster induced uncertainty (Shen et al., 2024), when combined with spatiotemporal features, has become an active research topic.

Overall, while existing studies have made substantial progress in revealing critical elements in urban transport networks, there remains significant scope for improvement in methods that simultaneously address large-scale network conditions, time-dependent characteristics, and efficient identification capability. This research context motivates further investigation into critical link identification approaches for time-dependent urban transport systems.

*1.3 Research objectives of this study*

In recent years, rapid advances in quantum computing and quantum optimisation have drawn increasing attention to their potential advantages in combinatorial optimisation and the analysis of large-scale complex systems (Ajagekar and You, 2019). In fields such as network design, network planning, and infrastructure optimisation, emerging studies have begun to introduce Quantum Annealing (QA) and Quantum Approximate Optimisation Algorithm (QAOA) frameworks into traditional network optimisation problems, to address the computational challenges posed by high dimensional decision spaces and complex constraints (Dixit and Niu, 2023; Dixit et al., 2024; Niu et al., 2025b). This development trend offers new technical possibilities for overcoming the efficiency limitations of classical methods in large-scale network analysis (Niu et al., 2025a).

Against this background, the objective of this study is to develop a critical link identification framework for large-scale urban transport networks that leverages the strengths of quantum optimisation in handling combinatorial complexity and captures the dynamic variation of link importance under time-dependent conditions. Specifically, this study aims to introduce a quantum modelling perspective to systematically identify network criticality across different traffic states, thereby revealing the key links that dominate operational risk in urban transport networks and the associated high risk time windows. The proposed framework seeks to provide new analytical tools and methodological support for resilience analysis and risk management of urban transport infrastructure.

*1.4 Research contributions*

This study investigates the identification of critical infrastructure elements in large-scale urban transport networks, with a focus on modelling methodology, empirical application, time-dependent characterisation, and computational feasibility. The main contributions are summarised as follows.

(1) From a modelling perspective, this study formulates the critical link identification problem in urban transport networks as an optimisation model that is suitable for solution on quantum hardware. This formulation provides a new modelling framework for characterising the importance of multi-link combinations under complex network operating conditions.

(2) From an empirical perspective, this study uses real world urban traffic data to systematically identify multiple critical link combinations that exert a dominant influence on network performance in a large-scale urban transport network. The results demonstrate the feasibility of applying the proposed framework at the city scale.

(3) From a time-dependent perspective, this study characterises the dynamic variation of critical link importance across the full time series, revealing the non-static nature of key network structures in urban transport systems over time.

(4) From a methodological and computational perspective, this study uses hybrid Constrained Quadratic Model (CQM) solver and quantum annealing on D-Wave hardware to achieve efficient solution of the critical link identification problem under large-scale network conditions and across the full temporal horizon. The results provide practical evidence of the applicability of this approach to urban transport infrastructure analysis.

Overall, the rest of this paper is organised as follows. **Section 2** develops a classical critical link identification model. **Section 3** reformulates the classical model into a quantum framework. **Section 4** presents an empirical analysis driven by real world urban network data. **Section 5** discusses the main findings and limitations of the study. **Section 6** concludes the paper.

**2. Classical model with time-dependent formulation**

This section formulates the classical critical link identification problem using real-world time-dependent traffic data. The model explicitly incorporates the temporal dimension through the time index $t$, while maintaining a single-period optimisation structure. The evaluation across all available time intervals is presented separately in the empirical section.

(1) Parameter definitions

The parameter settings and definitions used in this study are presented in **Table 1**.

**Table 1.** Table of parameter definitions in this study.

| | |
|---|---|
| $S$ | Set of road links, $S = \{1, 2, ..., s\}$ |
| $T$ | Set of observed time intervals, $T = \{1, 2, ..., t\}$ |
| $t_s(t)$ | Observed travel time on link $s$ at time interval $t$ |
| $Speed_s(t)$ | Observed speed on link $s$ at time interval $t$ |
| $Length_s$ | Physical length of link $s$ |
| $t_s^0$ | Free-flow travel time of link $s$ |
| $u_s$ | Binary decision variable: $u_s = 1$ if road link $s$ is disrupted, otherwise 0 |
| $t_s(t\|u_s)$ | Travel time of link $s$ under disruption state $u_s$ |
| $Delay_s(t\|u_s)$ | Delay on link $s$ under disruption state $u_s$ |
| $NDI(u, t)$ | Network delay index at time interval $t$ |
| $u$ | Decision vector, it collects the disruption decisions for all links in the network and serves as the decision variable, $u = (u_s)_{s \in S}$ |
| $k$ | Number of disrupted road links |

| $\gamma$ | Reflects the disruption severity (e.g., blockage, reduced capacity, or forced detours), $\gamma > 1$ |
|---|---|
| $c_s(t)$ | Single-link impact coefficient at time interval $t$ |
| $\beta_{sr}(t)$ | Interaction coefficient between links $s$ and $r$ at time interval $t$ |
| $\lambda$ | The penalty parameter enforcing the cardinality constraint, $\lambda > 0$ |

(2) Free-flow travel time

The free-flow travel time $t_s^0$ is determined from two independent perspectives:

Speed-based perspective:

$$t_s^{0,speed} = \frac{length_s}{v_s^{free}} \qquad (1)$$

where

$$v_s^{free} = \max_t speed_s(t) \qquad (2)$$

Minimum-time perspective:

$$t_s^{0,time} = \min_t t_s(t) \qquad (3)$$

Final definition:

$$t_s^0 = \min\left(t_s^{0,speed}, t_s^{0,time}\right) \qquad (4)$$

(3) Disrupted travel time

The travel time of link $s$ under disruption state $u_s$ at time interval $t$ is defined as:

$$t_s(t|u_s) = \begin{cases} t_s(t), & u_s = 0 \\ \gamma t_s(t) & u_s = 1 \end{cases} \qquad (5)$$

(4) Disrupted delay

$$Delay_s(t|u_s) = t_s(t|u_s) - t_s^0 \qquad (6)$$

This represents the time-dependent incremental congestion relative to free-flow conditions.

(5) Time-dependent network delay index

$$NDI(u,t) = \sum_{s \in S} Delay_s(t|u_s) \quad (7)$$

A larger value of $NDI(u,t)$ indicates higher congestion and greater vulnerability of the network at time interval $t$.

(6) Classical critical link identification problem

For each time interval $t$, the critical link identification problem is formulated as:

$$Z = \max_u NDI(u,t) \quad (8)$$

Subject to:

$$\sum_{s \in S} u_s = k \quad (9)$$

$$u_s = \{0,1\}, \forall s \in S, t \in T \quad (10)$$

This formulation seeks the set of links whose disruption produces the maximum system delay at time interval $t$. The model will subsequently be evaluated across all observed time intervals.

## 3. QUBO formulation of the time-dependent critical link identification model

The classical model identifies the set of links whose disruption maximises the time-dependent network delay index $NDI(u,t)$. To enable quantum optimisation, the constrained maximisation problem is reformulated as a Quadratic Unconstrained Binary Optimisation (QUBO) model for each time interval $t$. Following the general QUBO formulation principle in quantum optimisation (Glover et al., 2022; Dixit and Niu, 2023; Niu et al., 2025b), the QUBO formulation in this study is expressed as follows.

(1) General QUBO structure

The QUBO objective is written as:

$$Z_{qubo}(t) = \min_u H(u,t) \quad (11)$$

where the quadratic energy function is:

$$H(u,t) = -\sum_{s \in S} c_s(t) u_s - \sum_{\substack{s,r \in S \\ s \neq r}} \beta_{sr}(t) u_s u_r + \lambda \left( \sum_{s \in S} u_s - k \right)^2 \tag{12}$$

(2) Single-link impact coefficient $c_s(t)$

The single-link coefficient quantifies the marginal increase in system delay at time interval $t$ when only link $s$ is disrupted.

Let:

$NDI(u_s = 1, t)$: only link $s$ is disrupted at time interval $t$ (all other links remain normal);

$NDI(u_s = 0, t)$: all links are normal at time interval $t$.

Then

$$c_s(t) = NDI(u_s = 1, t) - NDI(u_s = 0, t) \tag{13}$$

(3) Interaction coefficient $\beta_{sr}(t)$

The interaction coefficient captures the nonlinear joint effect of disrupting links $s$ and $r$ at time interval $t$.

Let:

$NDI(u_s = 1, u_r = 1, t)$: links $s$ and $r$ are both disrupted, all other links normal.

$NDI(u_s = 1, u_r = 0, t)$: only link $s$ is disrupted.

$NDI(u_s = 0, u_r = 1, t)$: only link $r$ is disrupted.

$NDI(u_s = 0, u_r = 0, t)$: neither link $s$ nor link $r$ is disrupted.

Then

$$\beta_{sr}(t) = NDI(u_s = 1, u_r = 1, t) - NDI(u_s = 1, u_r = 0, t) - NDI(u_s = 0, u_r = 1, t) + NDI(u_s = 0, u_r = 0, t) \tag{14}$$

If $\beta_{sr}(t) > 0$: This indicates a negative synergistic effect, where the two failures exacerbate congestion.

If $\beta_{sr} < 0$: This indicates a degree of substitutability, meaning that the impact of one disruption partly offsets the other.

If $\beta_{st}(t) = 0$: The combined effect is exactly equal to the sum of the two independent effects, implying a purely additive relationship without interaction.

(4) Penalty coefficient

The penalty coefficient $\lambda$ must be sufficiently large so that any violation of the constraint produces a penalty dominating the objective value. A standard requirement is:

$$\lambda \gg \max\left(\sum_{s\in S}|c_s(t)|, \sum_{s,r\in S}|\beta_{sr}(t)|\right) \tag{15}$$

In practical computation, $\lambda$ is typically chosen to be one to two orders of magnitude larger than the maximum value among all QUBO coefficients, ensuring that constraint violations are always discouraged relative to objective contributions.

(5) Final QUBO model

For each time interval $t \in T$, the critical link identification problem is expressed in QUBO form as:

$$Z_{qubo}(t) = \min_{u} H(u,t) = -\sum_{s\in S}c_s(t)u_s - \sum_{\substack{s,r\in S \\ s\neq r}}\beta_{sr}(t)u_s u_r + \lambda\left(\sum_{s\in S}u_s - k\right)^2 \tag{16}$$

## 4. Numerical experiments with real-time traffic data

*4.1 Urban traffic data and network representation*

This study constructs a time-dependent transport network using real urban scale traffic operation records. The data cover the core urban area of the Greater Miami metropolitan area in southern Florida in the United States. The spatial boundary is defined by geographic coordinates, with latitude ranging from 25.1707 to 25.8714 and longitude ranging from −80.5067 to −80.1203. This region represents a typical coastal metropolitan transport corridor with a grid structured street network. The dataset used in this study contains a total of 1,048,575 link level observations (each represents the traffic condition of an individual road segment at a specific time step).

The original dataset is organised at the link observation level. The key attributes include origin and destination node identifiers, spatial location of each link expressed by latitude and longitude, travel time measured in minutes, speed measured in kilometres per hour, and a time step index. In terms of data completeness, no missing values are observed in the key fields. Regarding temporal coverage, the dataset is indexed by discrete time steps ranging from 0 to 665 at five-minute intervals (i.e. 0, 5, 10, …, 665). At each time step, a complete set of 13,704 link-level observations is recorded, representing a full snapshot of network conditions at that specific moment. Across the entire study period, this results in 134 distinct observation points, each corresponding to one time step and one complete network snapshot. In the subsequent analysis, the time step index is used to represent the underlying temporal progression, while the 134 distinct observation points constitute the effective temporal units for the time-dependent network analysis.

The observation period spans from 16:45:00 on 9 September 2017 to 04:00:00 on 10 September 2017, capturing the continuous evolution of traffic conditions from nighttime to the early hours of the following day. Network scale statistics indicate that the dataset used in this study comprises 6,081 nodes and 7,815 directed links. As shown in **Figure 1**, this study first visualises the spatial coverage of the network to illustrate its overall scale and distribution. Each point in the figure corresponds to the spatial sampling location of a link observation record, reflecting the spatial footprint of a large-scale urban transport network.

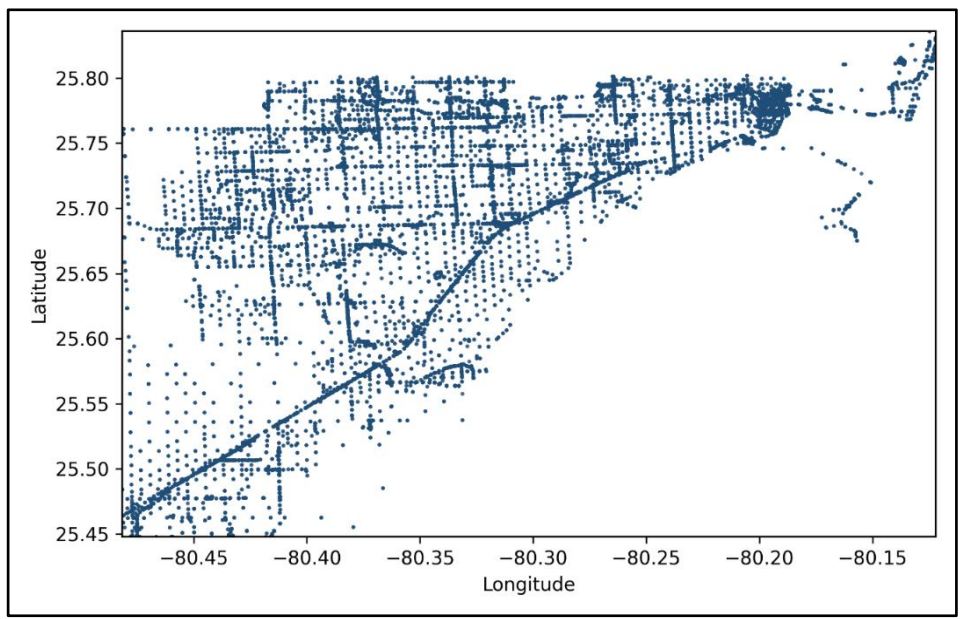

**Figure 1.** Spatial distribution of node latitude and longitude in the dataset used in this study.

Further, this study examines network connectivity to determine whether isolated sub networks or disconnected structures exist. The results show that the network is fully connected, meaning that all nodes belong to a single connected component, with no isolated nodes or independent sub networks. As shown in **Figure 2**, all nodes in the network fall within one giant connected component. This connectivity ensures the validity of the subsequent critical link identification analysis and avoids structural bias caused by data fragmentation.

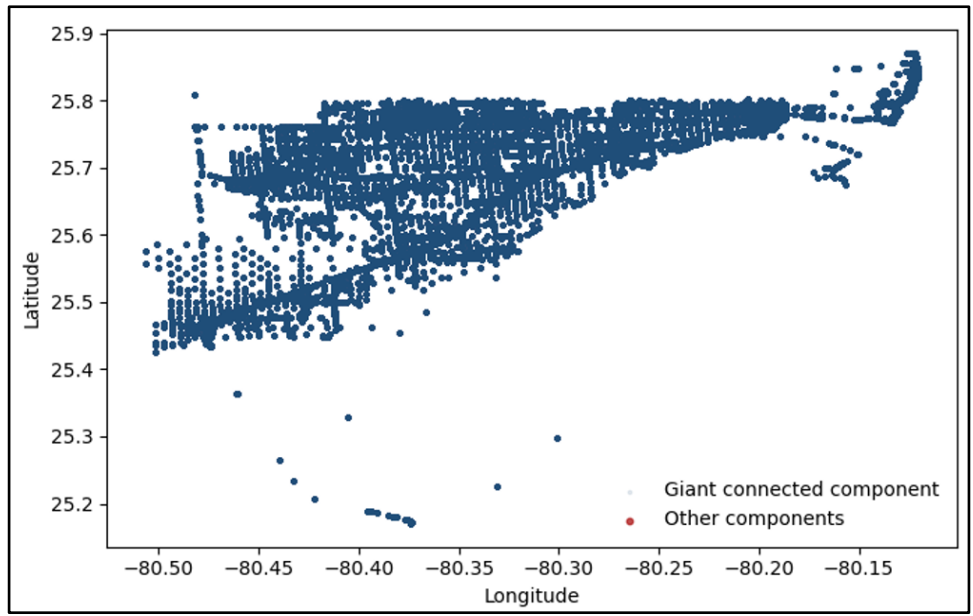

**Figure 2.** Network connectivity distribution under the real traffic dataset.

Overall, the data used in this study exhibit three key strengths. First, the sample size reaches the million level, which supports large-scale empirical analysis. Second, the temporal coverage is continuous with high resolution. Third, the network topology is complete and fully connected. This dataset and its network representation provide a reliable empirical foundation for subsequent QUBO oriented urban critical link identification and quantum optimisation experiments.

*4.2 Mapping observed traffic data to the QUBO model*

Based on the quantum QUBO model developed in this study, a critical next step is to map traffic observation data into the input parameters and energy function coefficients of the QUBO formulation. The overall mapping and computation workflow is illustrated in **Figure 3**. Link travel time, operating speed, and physical length recorded in the original data correspond respectively to the time-dependent link travel time, link speed, and link length in the model. Using these observed quantities, free flow travel time, perturbed travel time, and link delay are calculated step by step according to the model definitions and are then aggregated to form a time-dependent network delay metric.

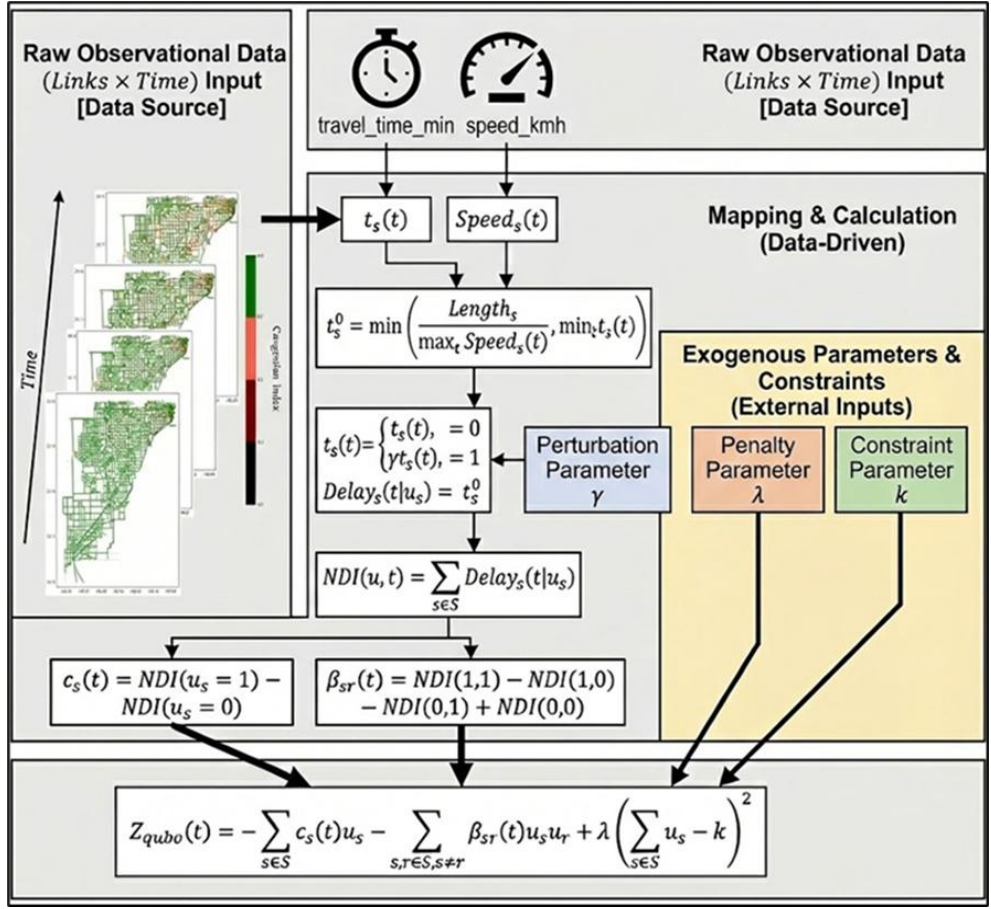

**Figure 3.** Data-driven mapping from observed traffic data to the QUBO model (This Figure was produced with assistance from *Gemini Pro*).

On this basis, the network delay index is evaluated under different disruption combinations to obtain single link impact coefficients and inter link interaction coefficients. These coefficients are then incorporated into the QUBO energy function at each corresponding time step, together with exogenously specified network delay index parameters, constraint scale, and penalty parameters. This procedure is executed consistently across the complete urban transport network and the full observation period, enabling the transformation and application of real traffic observation data within the QUBO optimisation model.

*4.3 Time-dependent critical link identification results and analysis*

*4.3.1 Overview of critical link identification results across the full time*

(1) Critical link identification results at a single time step

Using real traffic observation data from the regional road network in Florida in the United States, this study applies the hybrid CQM solver on the D-Wave hardware to solve the time-dependent critical link identification model formulated in QUBO form. In the model setting, constraints are imposed to control the number of links identified as critical at each time step, with the value fixed at $k = 20$.

Across all 665 discrete time steps, the model produces a total of 13,300-time step and critical link identification records, covering 7,814 directed links in the Florida regional road network. As shown in **Figure 4**, six consecutive time steps are randomly selected from the observation period, and the 20 critical links identified at each time step are visualised spatially. The results indicate that under a single time condition, the identified critical links exhibit a clear pattern of spatial concentration, with only a small number appearing in a dispersed manner. Most critical links are distributed along continuous major corridors and link segments connecting key nodes.

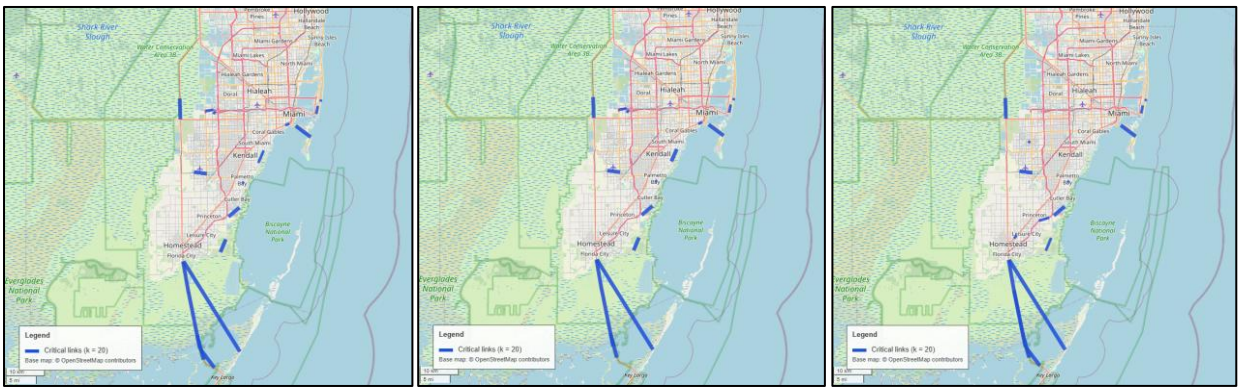

*Time_step = 90*   *Time_step = 180*   *Time_step = 270*

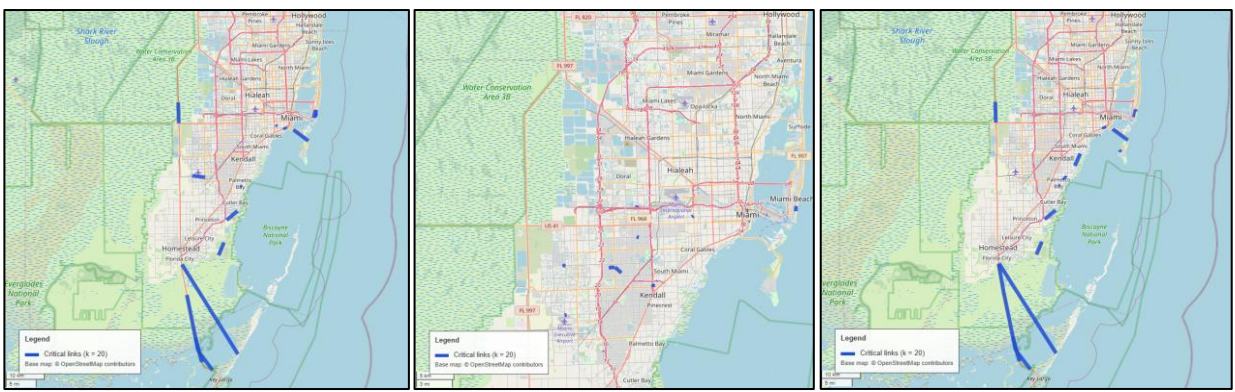

*Time_step = 360*   *Time_step = 450*   *Time_step = 540*

**Figure 4.** Spatial distribution of simultaneously disrupted critical links identified using the hybrid QUBO-CQM solution at selected time steps ($k=20$).

(2) Occurrence frequency characteristics of critical links

After obtaining critical link identification results for each time step, this study further computes the frequency with which each link is identified as critical over the entire observation period. The results show clear differences in temporal occurrence frequency across links. Some links are identified as critical only at a small number of time steps, while others repeatedly appear across many time steps, with frequencies substantially higher than those of most links in the network. As illustrated in **Figure 5**, high frequency critical links are mainly concentrated in highly urbanised areas of southern Florida, particularly within the Greater Miami metropolitan area and along major commuting corridors in its surrounding regions. These frequency statistics provide a basis for analysing the time-dependent spatial distribution of urban transport network vulnerability over a long-term horizon.

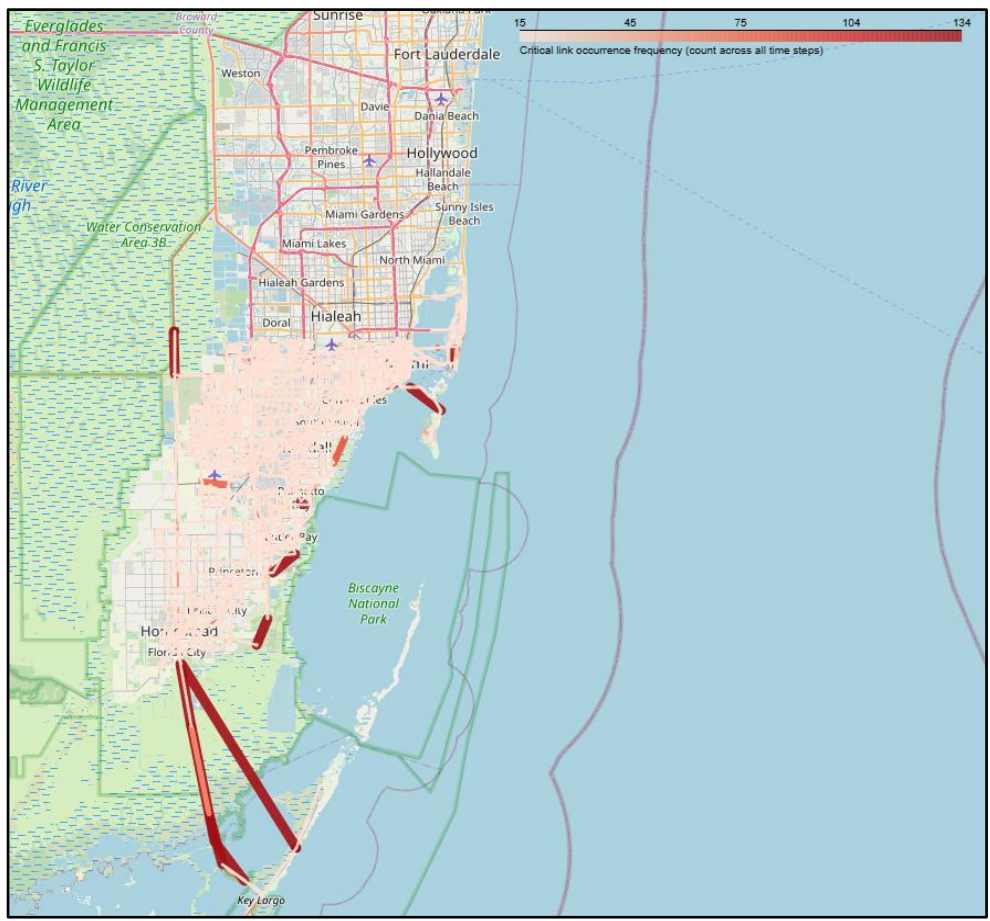

**Figure 5.** Occurrence frequency distribution of critical links over the full observation period.

By combining the critical link identification results at individual time steps with the temporal occurrence frequency characteristics of critical links, this study further selects the top 20 critical links with the highest occurrence frequencies across all 665 discrete time steps. **Table 2** reports the frequency characteristics of the selected critical link set, together with their spatial latitude and longitude coordinates and specific location coverage.

**Table 2.** Summary of critical link identification results over the observation period (Each critical link is represented by the geographic coordinates of its start and end nodes, where *Lat_from* and *Lon_from* denote the link origin and *Lat_to* and *Lon_to* denote the link destination in the directed network).

| ID | Frequency | Lat_from | Lon_from | Lat_to | Lon_to | location |
|----|-----------|----------|----------|--------|--------|----------|
| 1  | 134 | 25.23327588 | -80.43286386 | 25.18900246 | -80.3954724  | Monroe County, Florida, US |
| 2  | 134 | 25.36284318 | -80.46064742 | 25.4340965  | -80.47419505 | Miami-Dade County, Florida, US |
| 3  | 134 | 25.36385892 | -80.46100162 | 25.23327588 | -80.43286386 | Miami-Dade County, Florida, US |
| 4  | 134 | 25.20734379 | -80.4220178  | 25.26510112 | -80.43959148 | Monroe County, Florida, US |
| 5  | 134 | 25.80964286 | -80.48246737 | 25.7611379  | -80.48140955 | Miami-Dade County, Florida, US |
| 6  | 134 | 25.80976859 | -80.48225816 | 25.80976859 | -80.48225816 | Miami-Dade County, Florida, US |
| 7  | 134 | 25.22551634 | -80.33105758 | 25.22551634 | -80.33105758 | Monroe County, Florida, US |
| 8  | 134 | 25.22551634 | -80.33105758 | 25.43426468 | -80.47377399 | Miami-Dade County, Florida, US |
| 9  | 134 | 25.7922491  | -80.13343201 | 25.77452074 | -80.1374069  | Miami-Dade County, Florida, US |
| 10 | 134 | 25.48549629 | -80.36634611 | 25.45540994 | -80.37951376 | Miami-Dade County, Florida, US |
| 11 | 133 | 25.55793254 | -80.33140172 | 25.53732546 | -80.35895171 | Miami-Dade County, Florida, US |
| 12 | 133 | 25.74613639 | -80.1896234  | 25.72017685 | -80.15106794 | Miami-Dade County, Florida, US |
| 13 | 127 | 25.61826989 | -80.32477821 | 25.6146824  | -80.32047619 | Miami-Dade County, Florida, US |
| 14 | 112 | 25.75218692 | -80.20623507 | 25.74832758 | -80.21504427 | Miami-Dade County, Florida, US |
| 15 | 94  | 25.78709059 | -80.13170853 | 25.78709059 | -80.13170853 | Miami-Dade County, Florida, US |
| 16 | 86  | 25.69001612 | -80.26987579 | 25.66185587 | -80.28394494 | Miami-Dade County, Florida, US |
| 17 | 78  | 25.7911124  | -80.15498401 | 25.7911124  | -80.15498401 | Miami-Dade County, Florida, US |
| 18 | 77  | 25.63637467 | -80.4147552  | 25.64013125 | -80.4475288  | Miami-Dade County, Florida, US |
| 19 | 67  | 25.80761527 | -80.12483901 | 25.80255781 | -80.12648071 | Miami-Dade County, Florida, US |

| 20 | 63 | 25.78705643 | -80.13055212 | 25.79231428 | -80.13080997 | Miami-Dade County, Florida, US |

*4.3.2 Structural response of critical link combinations under varying disruption scales*

Building on the time-dependent critical link identification results presented above, this study further examines how critical link combinations respond to changes in disruption scale under a single traffic state. To this end, time step 475 is selected as a representative moment. While keeping the observed link travel times at this time step fixed, the number of simultaneously disrupted links $k$ is gradually increased to investigate the spatial structural evolution of the worst-case disruption scenario.

Specifically, at time step 475, this study sets the disruption scale $k = 20, 30, 40, 50, 60, 70, 80$, and solves the corresponding critical link combinations under different disruption scales using the same QUBO model structure and the hybrid solver on the D-Wave. It is important to emphasise that, in this analysis, only the disruption scale parameter $k$ is varied, while all other model parameters and input data remain unchanged.

As shown in **Figure 6**, when the disruption scale increases from $k = 30$ to $k = 50$, newly identified critical links do not replace the original ones. Instead, they extend progressively outward along existing critical corridors, forming larger and more continuous disruption structures. When $k = 60$, the spatial coverage of critical link combinations expands markedly, particularly along the Miami metropolitan area and the southern coastal corridors, where continuous multi-link disruption bands emerge. Notably, even when $k = 80$, the core links identified at earlier stages remain consistently selected. This finding suggests that under a fixed traffic state, the most sensitive critical link set exhibits a clear hierarchical structure, in which a small number of core links dominate worst-case simultaneous multi-link failure scenarios across different disruption scales, while additional links are incorporated gradually as the disruption scale expands.

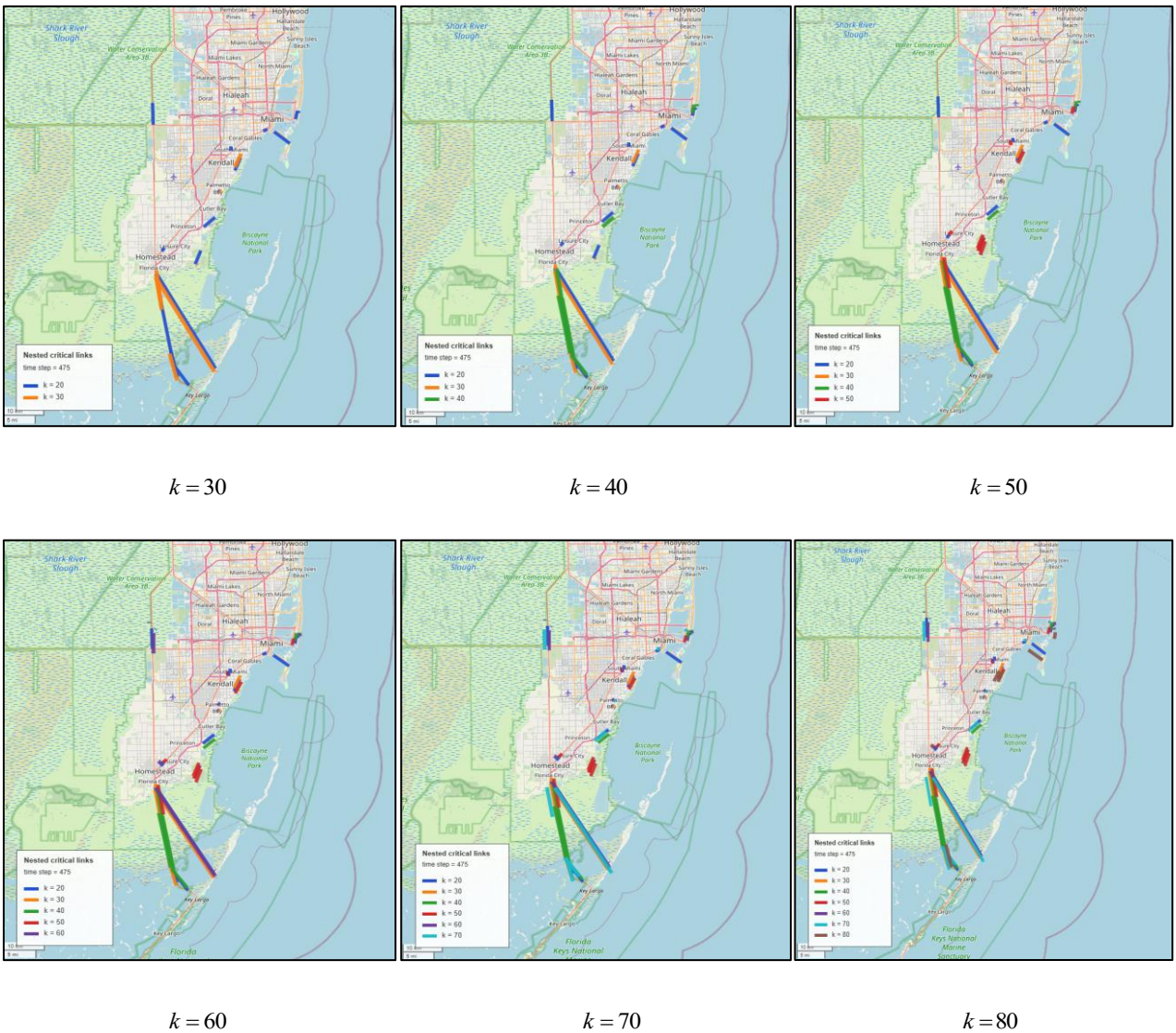

|  |  |  |
|---|---|---|
| $k = 30$ | $k = 40$ | $k = 50$ |
| $k = 60$ | $k = 70$ | $k = 80$ |

**Figure 6.** Spatial evolution of the network critical link identification structure.

*4.3.3 Time-dependent network delay index*

As shown in **Figure 7**, this figure presents the temporal evolution of the network delay index ( $NDI(u,t)$ ) under different disruption scales with $k = 20, 30, 40, 50, 60, 70, 80, 90, 100, 110$. The results indicate that as the disruption scale increases, the overall level of $NDI(u,t)$ rises markedly, demonstrating that an increase in the number of simultaneously failed links systematically amplifies the network delay index experienced by the network under the most adverse conditions.

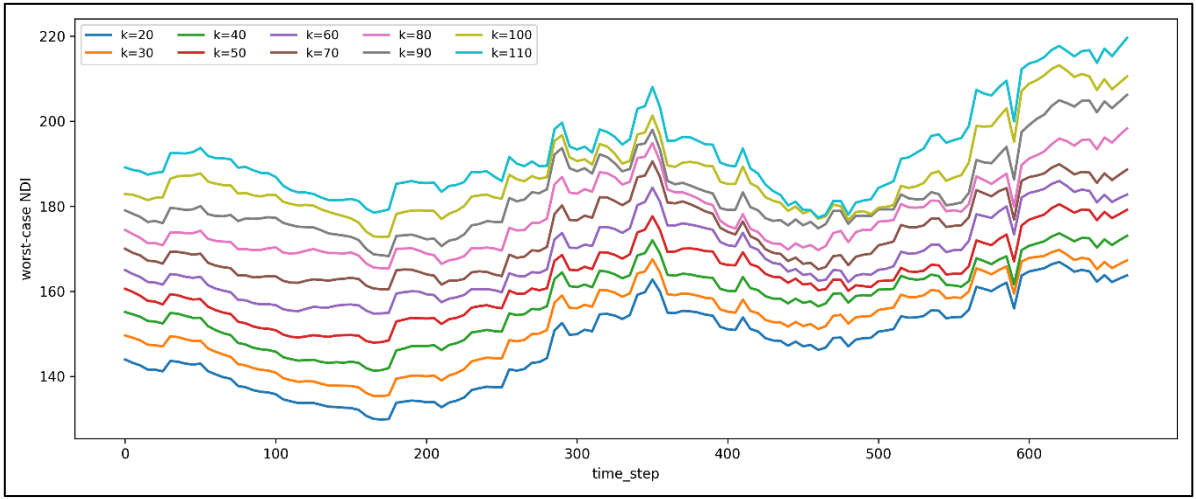

**Figure 7.** Time-dependent network delay index under different disruption scales.

Further inspection reveals that the *NDI*(*u*,*t*) curves associated with different disruption scales exhibit clear differences in the timing of peak values, the magnitude of fluctuations, and local variation patterns over time. This suggests that disruption scale not only affects the absolute level of network delay index but also alters the sensitivity of the network to disruptions under different traffic states. Overall, the results clearly show that the network delay index jointly constrained by traffic state and disruption scale, with a nonlinear coupling relationship between the two.

*4.3.4 Temporal risk patterns and high-risk time windows*

(1) Identification results of high-risk time windows

As shown in **Figure 8**, the identification results of high-risk time windows based on *NDI*(*u*,*t*) over time are presented, with the Top 10 percent and Top 5 percent high risk time steps highlighted. It can be observed that high risk time steps are not evenly distributed across the entire time series but are clearly concentrated in the later part of the observation period and tend to appear in consecutive sequences. This result indicates that the network faces periods of intensified disruption risk with strong temporal concentration. High risk states occur mainly within specific time windows rather than being randomly or uniformly distributed across the full observation horizon.

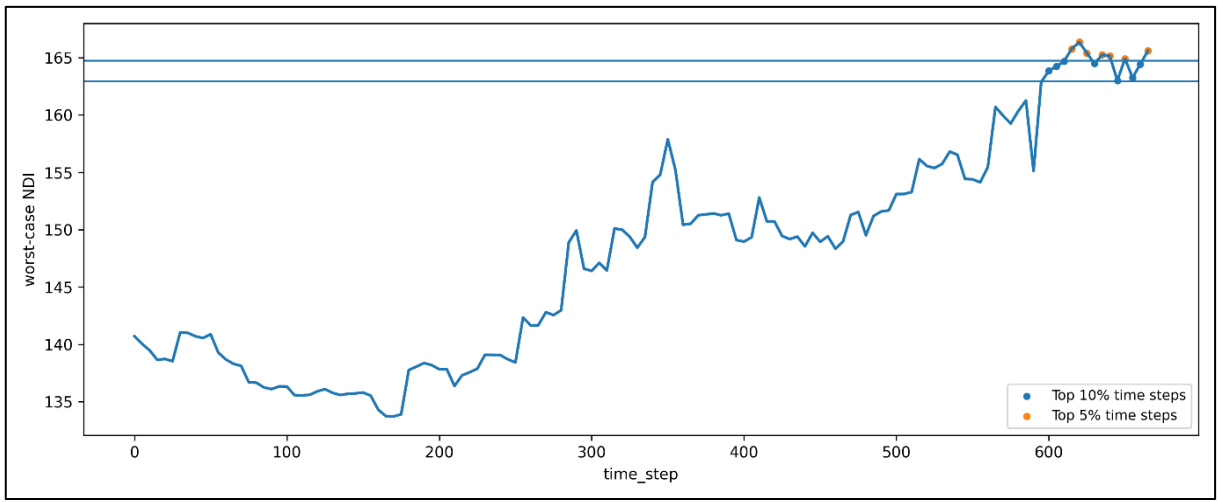

**Figure 8.** Identification of high-risk time windows based on $NDI(u,t)$ ( $k=20$ ).

(2) Analysis of change rate in network delay index

As shown in **Figure 9**, this study further presents the time series of the change rate of $NDI(u,t)$ between adjacent time steps, denoted as $\Delta NDI(u,t)$. For most time steps, $\Delta NDI(u,t)$ fluctuates slightly around zero, indicating that the network delay index generally remains relatively stable between consecutive traffic states. However, at a small number of time points, $\Delta NDI(u,t)$ exhibits pronounced positive or negative jumps, corresponding to rapid increases or decreases in $NDI(u,t)$.

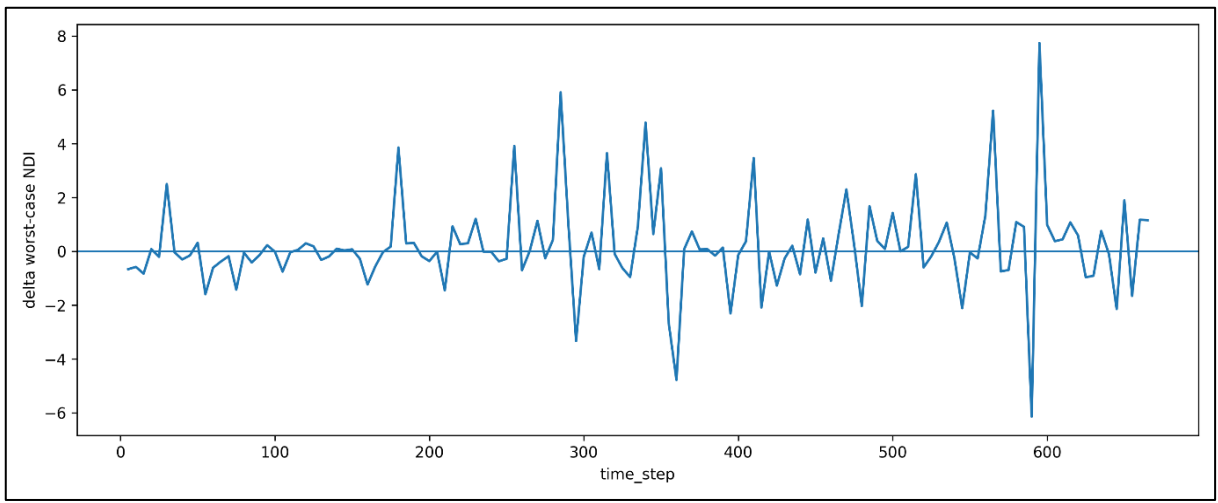

**Figure 9.** Temporal variation of changes in the worst-case Network Delay Index ( $k=20$ ).

When combined with the results shown in **Figure 8**, it becomes evident that these larger positive jumps tend to occur near high-risk time windows. This suggests that the transition into high-risk states is often accompanied by abrupt amplification of the network delay index.

5. Discussion

(1) Computation time and scalability analysis

As shown in **Figure 10**, the computation time required to evaluate the worst-case network delay index across different time steps exhibits a highly stable pattern overall. Apart from a small number of time steps, the computation time for a single time step remains within a narrow range of approximately 5 to 6 seconds, with no systematic growth trend as time progresses. This indicates that for most time steps, the proposed time-dependent worst-case disruption identification model maintains stable solution complexity, and its computational burden does not deteriorate noticeably with changes in network operating conditions.

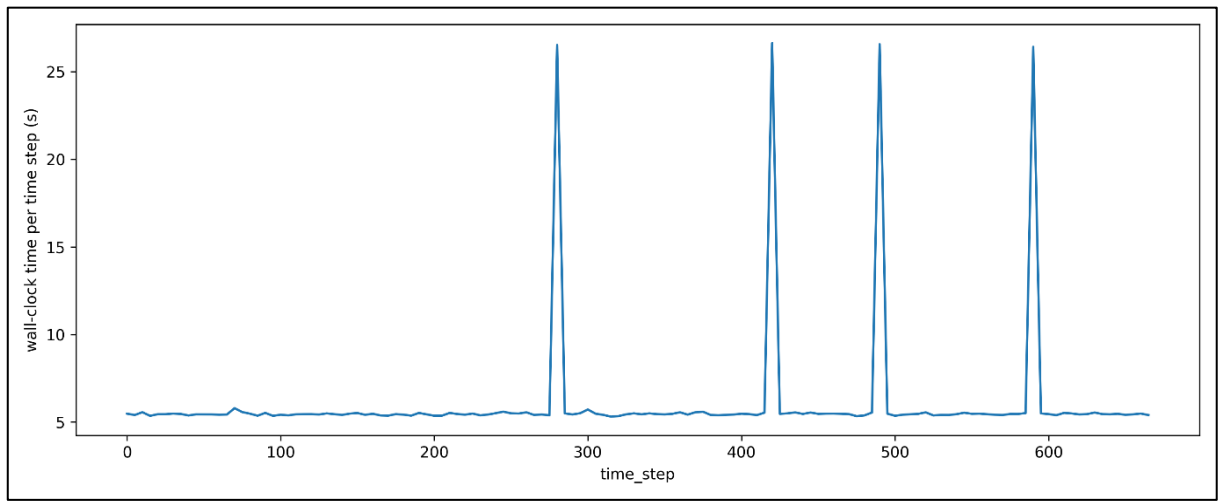

**Figure 10.** Solving time per time step using the hybrid CQM solver ( $k=20$ ).

(2) Interpretation of worst-case network delay index results

The vertical axes in **Figures 7** to **9** represent the worst-case network delay index, which directly corresponds to the objective value of the QUBO formulation, namely the minimum energy state. The results indicate that within high-risk time windows, even a limited scale of link failures can trigger a pronounced amplification of overall network delay, reflecting a substantial reduction in network redundancy and flow reallocation capability during these periods. By contrast, for most time steps, the worst-case delay level remains relatively stable and low, suggesting that under normal operating conditions the network retains a certain capacity to absorb disruptions. These findings imply that vulnerability assessments based solely on time-averaged indicators are likely to underestimate the true risk exposure during critical periods. Incorporating a time-dependent analysis of the worst-case delay therefore enables the identification of key time windows that dominate network performance degradation, providing a sound basis for time-sensitive risk management and targeted intervention strategies.

(3) Potential for real-time decision support

The second-level computation time achieved by the proposed framework suggests potential applicability in near real-time traffic analysis. By enabling rapid identification of critical links under different time-dependent network states, the framework can support timely understanding of structural risk in dynamic traffic conditions. While the method is not intended to directly generate control or routing decisions, it may provide structural guidance for traffic diversion, path re-planning, or incident response by highlighting links that require priority attention during high-risk periods.

As the disruption scale increases in **Figure 11**, the estimated computation time does not exhibit exponential growth. Instead, the solving time per time step increases gradually with $k$ and remains within the order of seconds even for large disruption scales. Occasional high computation peaks appear only at a small number of time steps and do not accumulate over time. This behaviour indicates that the proposed quantum-based framework maintains stable computational performance as problem size expands, supporting its applicability for real-time identification of critical links under dynamically evolving traffic conditions.

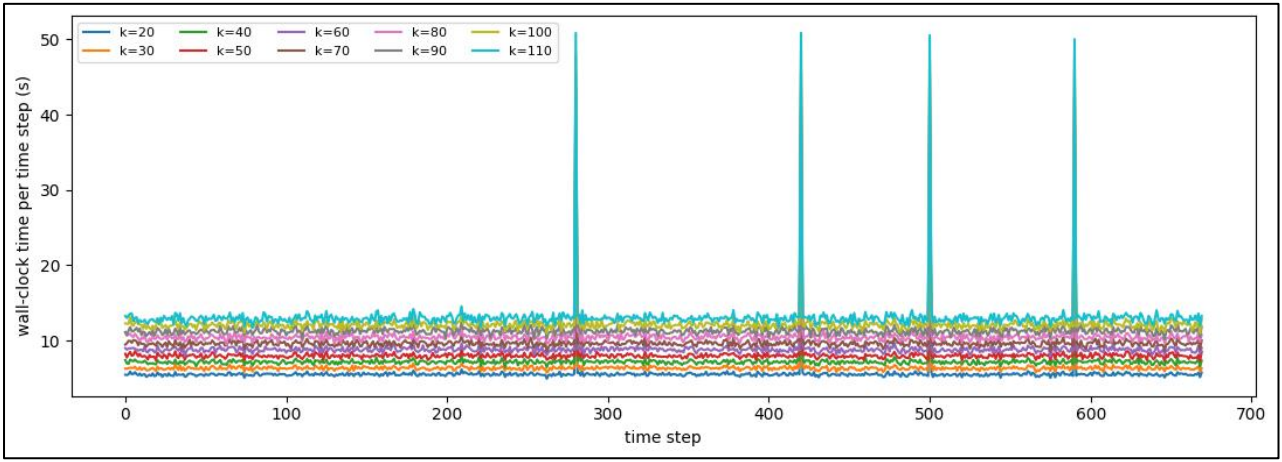

**Figure 11.** Solving time per time step for different disruption scales.

(4) Scope and limitations of the worst-case analysis

It should be noted that the worst-case network delay index examined in this study represents an upper-bound performance measure under constrained link disruption scenarios, rather than an estimate of expected or average network conditions. The results are therefore not intended to predict typical delay levels, but to identify time periods during which the network is structurally most vulnerable to adverse disruptions. In addition, the empirical analysis is conducted under a fixed disruption scale. While different values of this scale mainly affect the magnitude of the worst-case delay, the temporal patterns identified remain consistent. Direct comparison with traditional methods is not conducted because the problem setting addressed in this study differs fundamentally in both scale and formulation. Existing simulation-based or heuristic approaches are typically designed to evaluate average or expected network responses and become computationally prohibitive when extended to time-dependent, worst-case combinatorial disruption analysis on large-scale urban networks. Finally, the proposed framework is designed to complement, rather than replace, traditional traffic simulation approaches by providing a computationally tractable means of uncovering critical time windows and structural risk concentrations in large-scale urban transport networks.

## 6. Conclusion

This study proposes a time-dependent critical link identification framework for urban transport infrastructure. By reformulating the traditional critical link identification problem as a time-dependent QUBO optimisation model, the framework enables a systematic analysis of worst-case operating conditions in urban transport networks across the full temporal horizon. While retaining the logic of classical network modelling, the framework introduces a quantum annealing solution mechanism, which makes it possible to identify worst-case delay structures in large-scale and time-dependent networks.

Using high resolution traffic data from Florida in the United States, the empirical results demonstrate that urban transport network vulnerability exhibits strong time dependence. Critical links and their associated worst-case network delay index are not static. Instead, they vary substantially across different traffic states. High risk conditions concentrate within a small number of critical time windows. During these periods, even a limited scale of link failures can trigger a pronounced amplification of overall network delay. In contrast, during most normal operating time steps, the network maintains relatively stable and low delay levels even under worst-case disruption assumptions.

From a methodological perspective, the proposed time-dependent QUBO modelling framework highlights the potential of combining quantum optimisation techniques with real world urban traffic data. It provides an effective tool for revealing structural and temporal risks that traditional time averaged indicators struggle to capture. It is important to emphasise that the framework does not aim to predict routine traffic conditions. Instead, it functions as a stress testing and risk screening mechanism to identify critical infrastructure elements and critical time windows that dominate system performance degradation under extreme disruption scenarios. As quantum computing and hybrid optimisation technologies continue to advance, the proposed framework offers a scalable analytical foundation for urban scale infrastructure resilience assessment and the study of complex time-dependent systems.

**Author contributions**

**Junxiang Xu:** Writing – review & editing, Writing – original draft, Methodology, Formal analysis, Software, Conceptualization.

**Chence Niu:** Writing – review & editing, Methodology, Software.

**Tingting Zhang:** Writing – review & editing, Data processing.

**Divya Jayakumar Nair:** Review & editing, Supervision.

**Vinayak Dixit:** Review & editing, Supervision.

**Declaration of generative AI and AI-assisted technologies in the writing process**

During the preparation of this work, we employed ChatGPT-5.2 and Gemini to assist with language editing and polishing. No content was generated by AI. After using this tool, we thoroughly reviewed and edited the content as necessary and accept full responsibility for the content of the published article.